\documentclass[11pt]{article}
\usepackage{amsmath,amssymb,amsthm,enumerate,mathrsfs}
\usepackage[utf8]{inputenc}
\usepackage{empheq}
\usepackage{algpseudocode}
\usepackage{algorithmicx,algorithm}
\usepackage{mathtools}
\usepackage{color}
\usepackage{graphicx}
\usepackage{float}
\usepackage{caption}
\usepackage{tabularx,array}
\usepackage[nice]{nicefrac}
\usepackage[colorlinks=true, allcolors=blue]{hyperref}

\newtheorem{theorem}{Theorem}[section]
\newtheorem{lemma}[theorem]{Lemma}
\newtheorem{remark}[theorem]{Remark}

\newtheorem{definition}[theorem]{Definition}
\newcommand{\R}{\mathbb{R}}
\newcommand{\diag}{\mathrm{diag}}
\newcommand{\NESA}{{\sc Nesa}}
\newcommand{\NESAH}{{\sc Nesa$_{\tilde{H}}$}}

\newcommand{\boxedeq}[2]{\begin{empheq}[box={\fboxsep=6pt\fbox}]{align}\end{empheq}}
\algnewcommand{\Inputs}[1]{%
  \State \textbf{Inputs:}
  \Statex \hspace*{\algorithmicindent}\parbox[t]{.8\linewidth}{\raggedright #1}
}
\algnewcommand{\Initialize}[1]{%
  \State \textbf{Initialize:}
  \Statex \hspace*{\algorithmicindent}\parbox[t]{.8\linewidth}{\raggedright #1}
}
\algnewcommand{\Stop}[1]{%
  \State \textbf{Stopping Conditions:}
  \Statex \hspace*{\algorithmicindent}\parbox[t]{.8\linewidth}{\raggedright #1}
}

\newcommand{\MATLAB}{\textsc{Matlab}}

\numberwithin{equation}{section}

\title{Detecting negative eigenvalues of exact and approximate Hessian matrices 
in optimization}
\author{Warren Hare 
\thanks{Department of Mathematics, University of British Columbia, Okanagan 
Campus, Kelowna, B.C. V1V1V7, Canada 
(\texttt{warren.hare@ubc.ca}, ORCID 0000-0002-4240-3903).
Research partially supported by NSERC Discovery Grant \#2018-03865.  Research for this paper supported by France-Canada Research Funds 2022.
}
\and
Cl\'ement W. Royer
\thanks{
LAMSADE, CNRS, Universit\'e Paris Dauphine-PSL, 
Place du Mar\'echal de Lattre de Tassigny, 75016 Paris, France 
(\texttt{clement.royer@lamsade.dauphine.fr}, ORCID 0000-0003-2452-2172).
Research partially supported by Agence Nationale de la Recherche through 
program ANR-19-P3IA-0001 (PRAIRIE 3IA Institute).   Research for this paper supported by France-Canada Research Funds 2022.
}
}

\begin{document}

\maketitle

\begin{abstract}
Nonconvex minimization algorithms often benefit from the use of second-order 
information as represented by the Hessian matrix. When the Hessian  
at a critical point possesses negative eigenvalues, the corresponding eigenvectors can be used to search for further improvement in the objective function value. Computing such eigenpairs can be  computationally challenging, particularly if the Hessian matrix itself cannot 
be built directly but must rather be sampled or approximated. In blackbox 
optimization, such derivative approximations are built at a significant cost in terms of function values.

In this paper, we investigate practical approaches to detect negative 
eigenvalues in Hessian matrices without accessing the full matrix. We propose 
a general framework that begins with the diagonal and gradually builds 
submatrices to detect negative curvature. Crucially, our approach works both when exact Hessian coordinate values are available and when Hessian coordinate values are approximated.  
We compare several instances of our framework on 
a test set of Hessian matrices from a popular optimization library, and 
finite-differences approximations thereof. Our experiments highlight the 
importance of the variable order in the problem description, and show that 
forming submatrices is often an efficient approach to detect negative 
curvature.
\end{abstract}

\section{Introduction}
\label{sec:intro}

This paper considers unconstrained optimization problems of the form
\begin{equation}
\label{eq:pb}
	\min_{x \in \R^n} f(x),
\end{equation}
where $f:\R^n \rightarrow \R$ is a nonconvex $\mathcal{C}^2$ function ( $\mathcal{C}^2$ means twice  continuously differentiable). Due to the nonconvexity of the problem, it is desirable to design algorithms that converge towards second-order critical points, that is, points where the gradient vector is $0$ and the Hessian matrix does not possess negative curvature. Indeed, if the 
Hessian at a critical point has a negative eigenvalue, then better points (in the sense of having lower function values) can be found by moving along an 
eigenvector associated with this eigenvalue. This observation has long been 
a motivation for designing nonlinear optimization algorithms exploiting negative 
curvature~\cite{DGoldfarb_1980,GPMcCormick_1977,JJMore_DCSorensen_1979}, and 
recent advances have focused on schemes with worst-case complexity 
properties~\cite{CCartis_NIMGould_PhLToint_2012a,
FECurtis_ZLubberts_DPRobinson_2018,CWRoyer_SJWright_2018}.

Although the above methods possess strong theoretical guarantees, they are 
more expensive than other variants that do not rely on negative curvature 
exploitation. A standard way of detecting negative curvature consists of 
evaluating the full Hessian matrix and computing its minimum eigenvalue.  Both the evaluation and the eigenvalue calculation can be viewed as 
expensive procedures, especially in large dimensions. Matrix-free techniques 
for computing eigenvalues, such as Krylov subspace methods, have thus gained 
traction in the nonlinear optimization community, as they allow for computing
eigenvalue approximations based solely on Hessian-vector 
products~\cite{ARConn_NIMGould_PhLToint_2000,GFasano_SLucidi_2009}. This 
paradigm allows for more tractable variants of second-order methods, but 
still requires direct access to (partial) derivative information.

The situation becomes even more challenging while tackling 
problem~\eqref{eq:pb} using derivative-free (also called blackbox) optimization techniques, where 
the derivatives of the objective function are not directly  
employed~\cite{CAudet_WHare_2017,ARConn_KScheinberg_LNVicente_2009b}. In this setting, the Hessian matrix cannot be directly accessed and one must 
resort to approximations.  Hessian approximations can be constructed through second-order finite-difference formulas~\cite{JEDennis_RBSchnabel_1996}, although, some alternate estimation techniques have also been proposed~\cite{WHare_GJarryBolduc_CPlaniden_2020,WHare_KSrivastava_2020}. 
If the Hessian approximations are sufficiently accurate, then it is 
possible to design derivative-free techniques that exploit 
(approximate) negative curvature in order to converge to second-order 
stationary points~\cite{MAAbramson_LFrimannslund_TSteihaug_2014,
ARConn_KScheinberg_LNVicente_2009a,SGratton_CWRoyer_LNVicente_2016,
SGratton_CWRoyer_LNVicente_2020,DJudice_2015}. 

In derivative-free optimization, evaluating 
the objective function is often the computational bottleneck and building a full Hessian 
approximation requires a number of evaluations that scales quadratically 
with the dimension.  This scaling is not alleviated by the use of Krylov-type estimates, which scan the full matrix at every iteration. For this reason,  there remains a need for efficient 
negative curvature detection routines in a blackbox setting.

In this research, we explore numerical methods to rapidly determine if a 
negative eigenvalue exists.  We provide new tools that can be used to both check if second order optimality is (approximately) obtained (see Theorem \ref{thm:NESAHoutput}), and to rapidly determine a descent direction when negative eigenvalues exist (see Theorem \ref{thm:NESAHdescent}).  Our approach departs from standard numerical 
linear algebra techniques, such as Krylov subspace methods, in that it 
constructs submatrices one element at a time. As a result, this process is 
particularly well suited for the blackbox optimization setting, in which 
one can obtain an approximate Hessian coefficient at the expense of one or two 
additional function evaluations. We instantiate our framework based on 
several strategies for querying Hessian coefficients or their approximations,
which we validate on matrices extracted from the CUTEst optimization 
benchmark. Our experiments reveal that exploiting diagonal information to 
decide which submatrices to build can lead to faster detection of negative 
curvature. Perhaps surprisingly, we also provide numerical evidence that the 
order in which the variables are provided in optimization codes is often a 
good way of building submatrices, that allows for rapidly capturing significant 
negative curvature information. Our results illustrate the potential 
benefits of simple negative curvature estimates, even in a blackbox context.

The remainder of this paper is organized as follows. Section~\ref{ssec:eigs} 
concludes the introductory part of the paper by recalling some key results about eigenvalues and symmetric matrices. Section~\ref{sec:algo} describes our 
main algorithm and its variant tailored to derivative-free optimization. 
Numerical experiments with both variants are presented in 
Section~\ref{sec:exp}. Section~\ref{sec:conc} concludes the paper by 
discussing future uses of our approach.

\subsection{Background on eigenvalues of symmetric matrices}
\label{ssec:eigs}

Throughout this paper we only consider symmetric real-valued matrices, as our 
motivation stems from Hessian matrices in optimization over $\R^n$. Given a 
symmetric matrix $A \in \R^{n \times n}$, we say that $B$ is a principle 
submatrix of $A$ if $B$ can be constructed by deleting both the $i^{th}$ row 
and the $i^{th}$ column of $A$ for some values of $i$.  (Note that `some values of $i$' could consist of the empty set, as such $A$ is a principle submatrix of $A$.) Principle submatrices play a strong role in eigenvalue analysis, as is demonstrated by Cauchy's 
Interlacing Eigenvalue Theorem. This result is standard in the linear algebra 
literature, and has several extensions beyond the symmetric, real-valued 
setting~\cite[Theorem 10.1.1]{BNPartlett_1998}. It can also be applied when 
considering submatrices expressed in a different basis, as in the Lanczos' 
method~\cite[Chapter 13]{BNPartlett_1998}. In this paper, we will exploit 
the following form of the theorem.

\begin{theorem}[Cauchy's Interlacing Eigenvalue Theorem]\label{thm:cauchy}
Let $A \in \R^{n \times n}$ be a symmetric matrix and let 
$B \in \R^{m \times m}$ be a principal submatrix of $A$.  Suppose $A$ has 
eigenvalues $\lambda_1 \leq \lambda_2 \leq \ldots \leq \lambda_n$ and $B$ 
has eigenvalues $\beta_1 \leq \beta_2 \leq \ldots \leq \beta_m$.  Then, for 
each $k \in\{1, 2, \ldots, m\}$, we have 
\[
	\lambda_k \leq \beta_k \leq \lambda_{k+n-m}.
\]
\end{theorem}

An immediate corollary of Theorem \ref{thm:cauchy} is that a matrix is 
positive definite if and only if all principle submatrices are positive 
definite.\footnote{As $A$ is a principle submatrix of $A$, if all principle submatrices are positive definite, then obviously $A$ is positive definite.  Conversely, if $A$ is positive definite, then $\lambda_1 > 0$, so the smallest eigenvalue of any principle submatrix is strictly positive and thus all principle submatrices are positive definite.} In particular, if the $k^{th}$ eigenvalue of a principle submatrix is negative, then the $k^{th}$ eigenvalue of the original matrix is also negative.  

Consequently, if we seek to prove that negative eigenvalues exist, 
then it suffices to find a principle submatrix with a negative eigenvalue. 
This observation is the basis of our main algorithm, which is described in 
the next section.

\section{The Negative Eigenvalue Seeker Algorithm (\NESA)}
\label{sec:algo}

In this section, we present the Negative Eigenvalue Seeker Algorithm (\NESA), 
based on constructing principle submatrices in order to identify a submatrix 
with a negative eigenvalue, if any. Section~\ref{ssec:nesa} details the main 
version of the algorithm when exact matrices are used, while 
Section~\ref{ssec:nesah} is concerned with a variant of the algorithm 
dedicated to finite-difference matrix approximations. In both cases, we state 
the algorithm in a general fashion.  We describe several ways to instantiate 
the algorithms in Section~\ref{ssec:heuristics}.

\subsection{\NESA{} for exact Hessian matrices}
\label{ssec:nesa}

Algorithm~\ref{alg:NESA} (\NESA{}) describes our approach when the 
coefficients of $A$ can be evaluated directly. The idea behind \NESA{} is to 
gradually fill an auxiliary matrix $\tilde{A}$ with coefficients of $A$. This 
matrix is initialized to $\diag(A)$, i.e., 
the diagonal matrix with the same coefficients as $A$ on the diagonal and 
zeros elsewhere. If any element of the diagonal, which is a one-by-one submatrix, 
is sufficiently negative, then the method terminates immediately.
Otherwise, the algorithm selects a coordinate $(i,j)$ and updates $\tilde{A}$ 
by changing its $0$s in coordinates $(i,j)$ and $(j,i)$ to the value in those 
coordinates of $A$. All principle submatrices of $A$ that are stored in $\tilde{A}$ 
and include the coefficients $\tilde{A}_{i,j}$ and $\tilde{A}_{j,i}$ are then used 
to compute an estimate for the minimum eigenvalue of $A$. Depending on the value 
of that estimate, the process either repeats or stops.

Note that the stopping criterion involves a tolerance on the 
minimum eigenvalue of $A$. When this tolerance is set to $0$, as 
will be the case in our experiments, finding any nonpositive eigenvalue in a 
principle submatrix will lead to termination of the method.
%

\begin{algorithm}
\caption{Negative Eigenvalue Seeker Algorithm (\NESA)}\label{alg:NESA}
\begin{algorithmic}[1]
\Procedure{\NESA}{$A,\epsilon$}
\State \% $A \in \R^{n \times n}$ symmetric matrix that can be sampled 
coefficient-wise.
\State \% $\epsilon \ge 0$ stopping parameter (default $0$).
\State {\bf Initialize} $\tilde{A} = \diag(A)$ and 
$\lambda=\lambda_{\min}(\tilde{A})$.
\While{$\lambda \geq -\epsilon$ \& $\tilde{A} \neq A$}
    \State {\sc Select} $(i,j)$, $i>j$: a coordinate such that 
    $\tilde{A}_{i,j} \neq A_{i,j}$.  
    \State Update $\tilde{A}_{i,j} \leftarrow A_{i,j}$ and 
    $\tilde{A}_{j,i} \leftarrow A_{i,j}$.
    \State Compute the set $\mathcal{C}_{i,j}$ of all principle submatrices 
    of $A$ stored in $\tilde{A}$ that include $A_{i,j}$ and $A_{j,i}$.
    \State Set $\lambda = \min\{ \lambda_{\min}(C) | C \in \mathcal{C}_{i,j}\}$.
\EndWhile
\State {\bf Return} $\tilde{A}$ and $\lambda$.
\EndProcedure
\end{algorithmic}
\end{algorithm}

\begin{remark}
\label{rk:nesaimpl}
In practical implementation of lines 6 and 7, it is easiest to track which coefficients have been updated and which coefficients have not been updated; as opposed to selecting  coordinate such that $\tilde{A}_{i,j} \neq A_{i,j}$.  

Furthermore, we do not compute eigenvalues of all principle submatrices in 
$\mathcal{C}_{i,j}$, but restrict ourselves to the largest possible ones with 
no overlap. Finding such matrices amounts to finding the largest cliques in a 
graph, a well-understood problem in graph theory~\cite{CBron_JKerbosch_1973}.
\end{remark}

Regardless of the choice of $\epsilon$ and the way the indices are selected 
on Line 6 of Algorithm~\ref{alg:NESA}, we can provide guarantees regarding 
the termination and outputs of the method.  This is the purpose of the two 
following lemmas.

\begin{lemma}[Finite termination of \NESA{}]\label{lem:NESAterminates}
\NESA{} terminates after at most $n(n-1)/2$ iterations.
\end{lemma}

\proof Each iteration updates one coordinate $(i,j)$, $i>j$, with 
$\tilde{A}_{i,j} \neq A_{i,j}$.  As there are at most $n(n-1)/2$ such 
elements that are nonzero, we must have $\tilde{A} = A$ after at most 
$n(n-1)/2$ iterations.  
\qed

\begin{lemma}[Output of \NESA{}]\label{lem:NESAworks}
Upon termination of \NESA, the minimum eigenvalue of $A$ is not bigger than 
$\lambda$.
\end{lemma}

\proof If \NESA{} terminates due to $\tilde{A}=A$, then obviously 
$\lambda = \lambda_{\min}(A)$ and the statement is true. 
If \NESA{} terminates before $\tilde{A}=A$, the value of $\lambda$ is by 
definition the minimum eigenvalue of a principle submatrix $C$ of $A$. The 
result then follows directly from Theorem~\ref{thm:cauchy}.


\subsection{\NESA{} for approximate Hessians (\NESAH)}
\label{ssec:nesah}

As mentioned in the introduction, we are particularly interested in detecting 
negative curvature in a derivative-free setting, where approximate Hessian 
matrices are built using function values. In this paper, we focus on the most 
classical way of building such matrices via finite differences, as explained 
in the following definition.

\begin{definition}\label{def:DDHessian}
Let $f :\R^n \mapsto \R$ be $\mathcal{C}^2$ and $x \in \R^n$. For a given 
$h>0$, the finite-difference estimate of $\nabla^2 f(x)$ at $x$ is the matrix 
$\tilde{H}(x)$ such that its diagonal coefficients are given by
\begin{equation}
\label{eq:DDHessiandiag}
	\tilde{H}_{i,i}(x) = \frac{f(x+he_i)-2 f(x)+f(x-he_i)}{h^2} 
	\quad \mbox{ i=1,\dots,n,}
\end{equation}
and its off-diagonal coefficients are 
\begin{equation}
\label{eq:DDHessianoffdiag}
	\tilde{H}_{i,j}(x) = \frac{f(x+he_i+he_j)-f(x+he_i)-f(x+he_j)+f(x)}{h^2}
	\quad \mbox{ $1 \le j<i \le n$.}
\end{equation}
\end{definition}

Based on this approximation, we define a variant of \NESA{} dedicated to 
Hessian matrix approximations, called \NESAH{} (\NESA{} for approximate 
Hessians) and described in Algorithm~\ref{alg:NESAH}.

\begin{algorithm}
\caption{\NESA{} for approximate Hessians (\NESAH)}\label{alg:NESAH}
\begin{algorithmic}[1]
\Procedure{\NESAH}{$h, \epsilon$}
    \State \% $h>0$ finite difference parameter.
    \State \% $\epsilon \ge 0$ stopping parameter (default $0$).
	\State {\bf Initialize} set $\tilde{A}=\mathbf{0}$, then for $i=1, 2, \ldots, n$ set 
	$\tilde{A}_{i,i} = \tilde{H}_{i,i}(x)$ using~\eqref{eq:DDHessiandiag}. 
	Set $\lambda = \min\{\diag(\tilde{A})\}$.
	\While{$\lambda \geq -\epsilon$ \& 
	$\tilde{A} \neq \tilde{H}$}
    	\State {\sc Select} $(i,j)$, $i>j$: a coordinate such 
    	that $\tilde{A}_{i,j} \neq \tilde{H}_{i,j}$, where 
    	$\tilde{H}_{i,j}$ is given by~\eqref{eq:DDHessianoffdiag}.
    	\State Update $\tilde{A}_{i,j} \leftarrow \tilde{H}_{i,j}$ and 
    	$\tilde{A}_{j,i} \leftarrow \tilde{H}_{i,j}$.
    	\State Compute the set $\mathcal{C}_{i,j}$ of all principle submatrices 
    	of $\tilde{H}$ stored in $\tilde{A}$ that include $\tilde{H}_{i,j}$ and 
    	$\tilde{H}_{j,i}$.
    	\State Set $\lambda = \min\{ \lambda_{\min}(C) | C \in \mathcal{C}_{i,j}\}$.
	\EndWhile
	\State {\bf Return} $\tilde{A}$ and $\lambda$.
\EndProcedure
\end{algorithmic}
\end{algorithm}

\begin{remark}
\label{rk:nesahimpl}
Similar to \NESA, lines 6 and 7 of \NESAH{} are written using expressions 
checking $\tilde{A} \neq \tilde{H}$ and 
$\tilde{A}_{i,j} \neq  \tilde{H}_{i,j}$, while in practice, we simply keep track 
of which $(i,j)$ have been updated. Also similar to \NESA, it suffices to compute principle submatrices of maximum size only.
\end{remark}

Results analogous to those of the previous section can be established for 
\NESAH{}. Termination follows from the same argument than in the case of 
\NESA{}.

\begin{lemma}[Finite termination of \NESAH{}]\label{lem:NESAHterminates}
\NESAH{} terminates after at most $n(n-1)/2$ iterations.
\end{lemma}

On the other hand, interpreting the output of \NESAH{} requires more 
analysis than that of \NESA. First, we require a bound on the error between 
the eigenvalues of an approximated Hessian and the eigenvalues of the true 
Hessian.  In the next lemma, we use  $B_{h}(x)$  to denote the open ball of radius $h>0$ centred at $x$: $B_{h}(x) = \{ x' : \|x-x'\| < h\}$. 

\begin{lemma}[Error analysis for Hessian approximation]
\label{lem:HessianError}  
Let $f : \R^n \mapsto \R$ be $\mathcal{C}^2$. Let $x \in\R^n$ and $h>0$. 
Suppose that $\nabla^2 f$ is Lipschitz continuous on $B_{h}(x)$ with constant 
$L$.  If  $\tilde{H}$ is constructed as in Definition~\ref{def:DDHessian}, 
then
\begin{equation}
\label{eq:HessError}
	\|\nabla^2 f(x) - \tilde{H}\| \leq \frac{5}{3}\sqrt{n} L h
\end{equation}
and
\begin{equation}
\label{eq:EigError}
	|\lambda_{\min} (\nabla^2 f(x)) - \lambda_{\min}(\tilde{H})| 
	\leq  \frac{5}{3} \sqrt{n} L h.
\end{equation}
\end{lemma}

\proof 
Equation~\eqref{eq:HessError} is a standard result from numerical 
analysis~\cite[Lemma 4.2.3]{JEDennis_RBSchnabel_1996}. A bound on the error 
in approximating the Hessian directly leads to a bound on the error 
in approximating the minimum 
eigenvalue~\cite[Proposition 10.14]{ARConn_KScheinberg_LNVicente_2009b}, 
which immediately leads to equation \eqref{eq:HessError}.
\qed

We can now analyze the guarantees provided by the output of \NESAH.

\begin{theorem}[Output of \NESAH{}]
\label{thm:NESAHoutput}
Suppose that Algorithm~\ref{alg:NESAH} is applied to $f,x,h$ such that
$\nabla^2 f$ is Lipschitz continuous on $B_{h}(x)$ with constant $L$. 
Then, upon termination of \NESAH, the minimum eigenvalue of 
$\nabla^2 f(x)$ is not bigger than 
$\lambda + \frac{5}{3}\sqrt{n} L h$.
\end{theorem}

\proof By noting that \NESAH{} is simply \NESA{} applied to $\tilde{H}$, 
we see that Lemma \ref{lem:NESAworks} applies and the minimum eigenvalue of 
$\tilde{H}$ is not bigger than $\lambda$. Combining this 
with the error bound from Lemma~\ref{lem:HessianError} concludes the proof.
\qed

Theorem~\ref{thm:NESAHoutput} illustrates that detecting negative 
curvature in an approximate Hessian matrix can be due to the 
finite-difference approximation. In derivative-free optimization, the 
value of $h$ is often chosen in an adaptive fashion and tends to decrease 
as the algorithm unfolds: it is then possible to guarantee that ``true'' 
negative eigenvalues will be detected by the 
method~\cite{SGratton_CWRoyer_LNVicente_2016}.  

When the output of \NESAH{} is (sufficiently) negative, it can also be used to determine a descent direction from a first-order stationary point.

\begin{theorem}\label{thm:NESAHdescent}
Let $f : \R^n \mapsto \R$ be $\mathcal{C}^2$. Let $\bar{x} \in\R^n$ be a first-order stationary point of $f$ and $h>0$. Suppose that $\nabla^2 f$ is Lipschitz continuous on $B_{h}(\bar{x})$ with constant $L$.   Suppose \NESAH{} returns $\tilde{A}$ and $\lambda$ such that $\lambda < -\frac{5}{3} \sqrt{n} L h $.  Let $C$ be a principal submatrix of $\tilde{H}$ such that $\lambda = \lambda_{\min}(C)$ and $c$ be the  eigenvector of $C$ associated with $\lambda$.  Then $c$ defines a direction of decrease for $f$ at $\bar{x}$.
\end{theorem}

\proof Without loss of generality, we write
    $$\tilde{H} = \begin{bmatrix} C & X^\top \\ X & Z \end{bmatrix}$$
for some matrices $X$ and $Z$.  Define
    $$d = \begin{bmatrix}c \\ 0\end{bmatrix}.$$
By definition of $d$, we have
    $$d^\top \tilde{H} d = c^\top C c = \lambda \|c\|^2 = \lambda \|d\|^2.$$
Applying Taylor's theorem to $f$ at $\bar{x}$ and using that $\nabla f(\bar{x})=0$ by definition of $\bar{x}$, we have
    $$\begin{array}{rcl}
    f(\bar{x} + \tau d) &=& f(\bar{x}) + \tau \nabla f(\bar{x})^\top d + \frac{\tau^2}{2} d^\top \nabla^2 f(\bar{x}) d + O(\|\tau\|^3) \\
    &=&  f(\bar{x}) + \frac{\tau^2}{2} d^\top \left(\nabla^2 f(\bar{x})-\tilde{H}\right)d +\frac{\tau^2}{2} d^\top\tilde{H} d + O(\|\tau\|^3)\\
    &\leq&  f(\bar{x}) + \frac{\tau^2}{2} \|\nabla^2 f(\bar{x})-\tilde{H}\| \|d\|^2 +  \frac{\tau^2}{2}\lambda \|d\|^2 + O(\|\tau\|^3)\\
    &\leq& f(\bar{x}) + \frac{\tau^2}{2}\left(\frac{5}{3}\sqrt{n}L h + \lambda \right)\|d\|^2 + O(\|\tau\|^3).
    \end{array}$$
Applying $\lambda + \frac{5}{3} \sqrt{n} L h <0$ now shows that $f(\bar{x} + \tau d)  < f(\bar{x})$ for sufficiently small values of $\tau$.
\qed 

\subsection{Selecting coordinates}
\label{ssec:heuristics}

In this section, we describe several ways to instantiate both \NESA{} and 
\NESAH{} by defining strategies to select coordinates in the process of 
building submatrices. 

Given an $n \times n$ matrix, we define such a selection strategy using 
two ingredients. The first ingredient is a permutation of the set 
$\{1, 2, \ldots, n\}$, while the second ingredient describes how this 
permutation is used to define coordinates to select.

We first describe the second ingredient of our approach. Given a 
permutation $P=[p_1, p_2, \ldots, p_n]$ of $\{1,\dots,n\}$, we 
design two different ways of creating a selection order called 
Build 1 and Build 2, respectively described by Algorithm~\ref{alg:build1} 
and Algorithm~\ref{alg:build2}.  The output of these algorithms ($Order$) is an ordered list of coordinates that is used to determine which $(i,j)$ to select in line 6 of Algorithm \ref{alg:NESA} (respectively line 8 of Algorithm \ref{alg:NESAH}).

\begin{algorithm}[ht!]
\caption{Build 1}\label{alg:build1}
\begin{algorithmic}[1]
\Procedure{Build1}{$P$}
    \State \% $P$ a permutation of the set $\{1, 2, \ldots, n\}$
    \State $\begin{array}{rl} 
    Order = & 
    [(p_1,p_2), (p_1,p_3), (p_1,p_4), \ldots, (p_1, p_n),\\
    & \quad (p_2, p_3), (p_2, p_4), \ldots, (p_2,p_n), \ldots, (p_{n-1},p_n)]
    \end{array}$
    \State {\bf Return} $Order$
\EndProcedure
\end{algorithmic}
\end{algorithm}

\begin{algorithm}[ht!]
\caption{Build 2}\label{alg:build2}
\begin{algorithmic}[1]
\Procedure{Build2}{$P$}
    \State \% $P$ a permutation of the set $\{1, 2, \ldots, n\}$
    \State $\begin{array}{rl} Order = & 
    [(p_2,p_1), (p_3,p_2), (p_3,p_1), (p_4, p_3) \ldots, (p_4, p_1), \ldots, \\
    & \quad (p_n, p_{n-1}), (p_n, p_{n-2}), \ldots, (p_n,p_1)]
    \end{array}$
    \State {\bf Return} $Order$
\EndProcedure
\end{algorithmic}
\end{algorithm}

To understand the difference between these two algorithms, suppose that 
we apply \NESA{} to a $4 \times 4$ matrix using $P=\{1,2,3,4\}$. As shown in 
Figure~\ref{fig:build1example}, Build 1 focuses on one row at a time and 
expands that row until the row is complete.  It then moves to the next available 
row and repeats the process.  Effectively, this is creating many small principle 
submatrices at the beginning of the process and then slowly merging them as the 
process continues.
On the other hand, as illustrated in Figure~\ref{fig:build2example}, Build 2 
starts near the diagonal and builds outwards until it the row is complete.  
Effectively, this creates a principle submatrix and then focuses on expanding the 
size of that principle submatrix as rapidly as possible.

\begin{figure}[ht!]
$$\begin{bmatrix}
x_1 & & & \\
& x_2 & & \\
& & x_3 & \\
& & & x_4
\end{bmatrix}
\longrightarrow
\begin{bmatrix}
x_1 & 1 & & \\
1 & x_2 & & \\
& & x_3 & \\
& & & x_4
\end{bmatrix}
\longrightarrow
\begin{bmatrix}
x_1 & 1 & 2 & \\
1 & x_2 & & \\
2 & & x_3 &  \\
& & & x_4
\end{bmatrix}
\longrightarrow
\begin{bmatrix}
x_1 & 1 & 2 & 3 \\
1 & x_2 &  & \\
2 & & x_3 &  \\
3 & & & x_4
\end{bmatrix}$$
$$\longrightarrow
\begin{bmatrix}
x_1 & 1 & 2 & 3\\
1 & x_2 & 4 & \\
2 & 4 & x_3 & \\
3 & & & x_4
\end{bmatrix}
\longrightarrow
\begin{bmatrix}
x_1 & 1 & 2 & 3\\
1 & x_2 & 4 & 5 \\
2 & 4 & x_3 & \\
3 & 5 & & x_4
\end{bmatrix}
\longrightarrow
\begin{bmatrix}
x_1 & 1 & 2 & 3\\
1 & x_2 & 4 & 5\\
2 & 4 & x_3 & 6\\
3 & 5 & 6 & x_4
\end{bmatrix}
$$\caption{An example of Build 1 working with \NESA{} will fill out a 
$4 \times 4$ matrix with $P=\{1, 2, 3,4\}$.  The $x = [x_1, x_2, x_3, x_4]^\top$ on the diagonal are filled out during initialization. Each number represents the iteration of 
\NESA{} at which the corresponding coefficient would be updated.}\label{fig:build1example}
\end{figure}

\begin{figure}[ht!]
$$\begin{bmatrix}
x_1 & & & \\
& x_2 & & \\
& & x_3 & \\
& & & x_4
\end{bmatrix}
\longrightarrow
\begin{bmatrix}
x_1 & 1 & & \\
1 & x_2 & & \\
& & x_3 &  \\
& & & x_4
\end{bmatrix}
\longrightarrow
\begin{bmatrix}
x_1 & 1 &  & \\
1 & x_2 & 2 & \\
 & 2 & x_3 &  \\
 & & & x_4
\end{bmatrix}
\longrightarrow
\begin{bmatrix}
x_1 & 1 & 3 & \\
1 & x_2 & 2 & \\
3 & 2 & x_3 &  \\
 & & & x_4
\end{bmatrix}$$
$$\longrightarrow
\begin{bmatrix}
x_1 & 1 & 3 & \\
1 & x_2 & 2 & \\
3 & 2 & x_3 & 4 \\
 & & 4 & x_4
\end{bmatrix}
\longrightarrow
\begin{bmatrix}
x_1 & 1 & 3 & \\
1 & x_2 & 2 & 5 \\
3 & 2 & x_3 & 4 \\
 & 5 & 4 & x_4
\end{bmatrix}
\longrightarrow
\begin{bmatrix}
x_1 & 1 & 3 & 6\\
1 & x_2 & 2 & 5\\
3 & 2 & x_3 & 4\\
6 & 5 & 4 & x_4
\end{bmatrix}
$$\caption{An example of how Build 2 working with \NESA{} fills out a 
$4 \times 4$ matrix with $P=\{1, 2, \ldots, n\}$.  The $x=[x_1, x_2, x_3, x_4]^\top$ on the diagonal are 
filled out during initialization. Each number represents the iteration of 
\NESA{} at which the corresponding coefficient would be updated.}\label{fig:build2example}
\end{figure}

We now elaborate on the first ingredient of our strategies, namely the choice 
of a permutation of $\{1,\dots,n\}$, where $n$ is the problem dimension. 
Since our algorithm begins with the matrix diagonal, we seek to exploit this 
information upon selecting the submatrices to be formed (recall that we only 
form submatrices when all diagonal elements are nonnegative). A naive approach 
consists in using the natural order of the indices. Following previous 
strategies in derivative-free 
optimization~\cite{SGratton_CWRoyer_LNVicente_2016}, another possibility is 
to consider the smallest diagonal elements as more promising for building 
submatrices with negative curvature, and to give priority to the associated 
indices. Conversely, one could prioritize the largest
coefficients, in the hopes that changes regarding these coefficients might have the largest impact.  Finally, 
combining indices corresponding to the smallest and the largest coefficients 
could also be beneficial. These considerations lead us to the following four 
heuristics for the permutation $P$:
\begin{enumerate}
	\item {\tt Ordered}: $P=[1, 2, 3, \ldots, n]$;
	\item Smallest to Largest Diagonal Element ({\tt S2Lde}): 
	Choose $P$ such that 
	\[ A_{p_1, p_1} \leq  A_{p_2, p_2} \leq \ldots \leq A_{p_n, p_n}.\]
	\item Largest to Smallest Diagonal Element ({\tt L2Sde}): 
	Choose $P$ such that 
	\[ A_{p_1, p_1} \geq  A_{p_2, p_2} \geq \ldots \geq A_{p_n, p_n}.\]
	\item Interlacing Diagonal Elements ({\tt Ide}): Create $P^{\tt temp}$ 
	as in {\tt S2Lde}, then set 
	\[ 
		P = [p^{\tt temp}_1, p^{\tt temp}_n, p^{\tt temp}_2, 
		p^{\tt temp}_{n-1}, \ldots, p^{\tt temp}_{\lceil n/2 \rceil}],
	\]
	where $\lceil \cdot \rceil$ rounds-up to the nearest integer.
\end{enumerate}

Note that the {\tt Ordered} strategy is the only one that does not 
leverage information from the matrix. In fact, if the matrix $A$ is 
randomly generated, there is no reason to believe that the order of 
the coefficients matters. However, as we will see in the next section, 
this turns out not to be the case on matrices coming from optimization 
benchmarks.

\section{Experiments using exact and approximate Hessian matrices}
\label{sec:exp}

In this section, we investigate the numerical behavior of \NESA{} and 
\NESAH{}, equipped with the heuristics described in 
Section~\ref{ssec:heuristics}. We compare those heuristics on a 
matrix test set formed using Hessian matrices arising in the CUTEst 
library~\cite{NIMGould_DOrban_PhLToint_2015}.  All implementations are in \MATLAB (version 9.10.0.1602886, R2021a), and eigenvalue computations are done using the {\tt eig.m} command.  All implementations are available on github\footnote{\href{https://github.com/clementwroyer/negative-eigs}
{https://github.com/clementwroyer/negative-eigs}}.

\subsection{Test Problems}
\label{ssec:pbms}

To construct our test problems, we begin with a subset of 49 unconstrained 
problems from the CUTEst collection~\cite{NIMGould_DOrban_PhLToint_2015} for 
which the objective function is twice continuously differentiable and the 
Hessian matrix at the initial point has a negative eigenvalue. The complete 
problem list, along with their dimensions and classifications, is given in 
Table~\ref{tab:listpbms}. 

Despite the diverse nature of the test set, we point out that all of these 
problems were created by a human. As we will see later, this is relevant to the performance of our method.

{\small
\begin{table}[ht!]
\caption{List of the \texttt{CUTEst} test problems.}
\begin{center}
\begin{tabular}{|c|c|c||c|c|c|}
	\hline
	Name &Dimension &Nature &Name &Dimension &Nature\\
	\hline
	ALLINITU    &4  	&Academic 	&BIGGS6		&6  	&Academic \\
	BOX3        &3 		&Academic	&BRYBND		&10		&Academic \\
	DENSCHND    &3		&Academic	&DENSCHNE	&3		&Academic\\
	DIXMAANA    &15 	&Academic	&DIXMAANB	&15		&Academic\\
	DIXMAANC    &15		&Academic 	&DIXMAAND	&15		&Academic\\
	DIXMAANE    &15		&Academic 	&DIXMAANF	&15		&Academic\\
	DIXMAANG    &15		&Academic 	&DIXMAANH	&15		&Academic\\
	DIXMAANI    &15		&Academic 	&DIXMAANJ	&15		&Academic\\
	DIXMAANK    &15		&Academic 	&DIXMAANL	&15		&Academic\\
	ENGVAL2     &3		&Academic	&EXPFIT		&2		&Academic\\
	FMINSURF    &16  	&Modeling	&FREUROTH	&10		&Modeling\\
	GROWTHLS    &3  	&Academic 	&GULF		&3		&Modeling\\
	HAIRY       &2		&Academic	&HATFLDD	&3		&Academic\\
	HATFLDE     &3		&Academic	&HEART6LS	&6		&Modeling\\
	HEART8LS    &8		&Modeling	&HELIX		&3		&Academic\\
	HIMMELBB    &2		&Academic	&HIMMELBG	&2		&Academic\\
	HUMPS       &2		&Academic	&KOWOSB		&4		&Modeling\\
	LOGHAIRY    &2		&Academic	&MEYER3		&3		&Real\\
	MSQRTALS    &4		&Academic	&MSQRTBLS	&9		&Academic\\
	OSBORNEA    &5		&Modeling	&OSBORNEB	&11		&Modeling\\
	PENALTY3    &50		&Academic	&SCOSINE	&10		&Academic\\ 
	SINQUAD     &50		&Academic	&SPARSINE	&10		&Academic\\ 
	SPMSRTLS    &28		&Academic	&VAREIGVL	&10		&Academic\\
	VIBRBEAM	&8		&Modeling	&WATSON		&12		&Academic\\
	YFITU       &3		&Modeling	& & & \\
	\hline
\end{tabular}
\end{center}
\label{tab:listpbms}
\caption{Problem list with CUTEst 
classification~\cite{NIMGould_DOrban_PhLToint_2015}.
Academic nature means that the problem was constructed by researchers for the academic purpose of testing one or more algorithms.  Modeling nature means that the problems was constructed as part of a modeling exercise, but the solution is not used in a genuine practical application.  Real nature means that the problem was constructed from an application for purposes other than testing algorithms.
}
\end{table}
}

For each test problem, we ran two iterations of Newton's method to generate 
points $x^0$, $x^1$, and $x^2$, where $x^0$ is the initial point provided in 
the CUTEst collection. This expands our test set to 147 matrices, of which 134 
possess negative curvature. Out of these 134 matrices, we discard 52 matrices 
that have negative diagonal elements, since in that case both \NESA{} and 
\NESAH{} terminate without iterating. 

\subsection{Using \NESA{} on exact Hessian matrices}
\label{ssec:numnesa}

In our first experiment, we apply \NESA{} to the 82 Hessian matrices 
computed using the procedure described in the previous paragraph with 
$\epsilon=0$. We compare Build 1 and Build 2 using the 4 permutation heuristics for 
\NESA{}, leading us to 8 different methods. For each problem, we determine 
which method(s) used the least number of iterations to detect a negative 
eigenvalue. Table~\ref{tab:exactALL} summarizes the results: we notice that the
percentages are similar for both build types. In both cases, we see a sharp 
advantage to selecting coordinates in the standard order $\{1,\dots,n\}$ as 
well as using the smallest-to-largest diagonal elements ({\tt S2Lde}). 

\begin{table}[ht!]
\begin{center}
\begin{tabular}{|l|c|c|c|c|c|}\hline
	 & {\tt Ordered} & {\tt S2Lde} & {\tt L2Sde} & {\tt Ide} \\
	\hline
	Build 1 &59.8 &50.0 &26.8 &35.4 \\
	Build 2 &58.5 &56.1 &31.7 &41.5 \\
    \hline
\end{tabular}
\end{center}
\caption{Percentage of problems where each build type and permutation 
heuristic resulted in the least number of iterations to detect a negative 
eigenvalue (82 matrices).}
\label{tab:exactALL}
\end{table}

When the problem dimension is small, the number of possible strategies 
actually exceeds the number of possible orderings that can be used. Indeed, 
if $n=2$, then there is actually only 1 possible order for the selection 
procedure.  There are $3!=6$ possible orders for the selection procedure when 
$n=3$.  (In general, when the dimension is $n$, then there are 
$\left[(n-1)n/2\right]!$ possible orders for the selection procedure.)  
Therefore, we also present the results after discarding two- and three- 
dimensional problems (leaving 63 matrices out of 85). The updated percentages 
are given in Table~\ref{tab:exact4up} and we see that Heuristic 1, {\tt Ordered}, again results in the overall best performance. 

\begin{table}[ht!]
\begin{center}
\begin{tabular}{|l|c|c|c|c|c|}\hline
	 & {\tt Ordered}: & {\tt S2Lde} & {\tt L2Sde} & {\tt Ide} \\
	\hline
	Build 1 &50.0 &38.3 &13.3 &13.3 \\
	Build 2 &48.3 &46.7 &20.0 &21.7 \\
    \hline
\end{tabular}
\end{center}
\caption{Percentage of problems of dimension 4 and higher where each build type 
and permutation heuristic resulted in the least number of iterations to 
detect a negative eigenvalue (60 matrices).}
\label{tab:exact4up}
\end{table}

Before analyzing Table \ref{tab:exact4up}, we remark that the best variant of \NESA{} detects negative curvature within $2$ iterations for 57 out of 82 matrices, which is significantly faster than the upper bound provided in Lemma \ref{lem:NESAterminates}. The worst performance is observed on a matrix from 
problem VAREIGVL (dimension 10): \NESA{} terminates in $28$ iterations, which still improves over the 
theoretical maximum of $(10\times 9)/2 = 45$ iterations. 

Examining Table \ref{tab:exact4up}, perhaps the most surprising result is the performance of the {\tt Ordered} heuristic, particularly since the {\tt Ordered} heuristic does not appear to use any information about the problem to make its decision. However, we conjecture that the {\tt Ordered} heuristic is in fact using a very powerful piece of information about the 
problem: when creating test problems, researchers naturally order the 
variables from {\em most impactful} to {\em least impactful}.  We further 
argue that this is even more pronounced when practitioners are constructing 
real-world problems, as models are naturally built starting from the {\em 
most impactful} variables.
To check our hypothesis, we thus consider the test set of 60 matrices, but 
we apply a random permutation matrix to all matrices in order to re-order the 
variables prior to applying our heuristics. The results appear in 
Table~\ref{tab:exactperm}, and show that the ordering $\{1,\dots,n\}$ is no 
longer beneficial compared to our heuristic consisting in selecting the 
lowest-curvature coordinate. We note, however, that the percentages are slightly in favor of 
Build 2, especially for the {\tt Ordered} and {\tt S2Lde} heuristics. Recall that 
Build 2 promotes the construction of large submatrices, that are likely to be used by 
our algorithm for better eigenvalue estimation. The results of Table~\ref{tab:exactperm} 
therefore suggest that Build 2 is a better strategy on average.

\begin{table}[ht!]
\begin{center}
\begin{tabular}{|l|c|c|c|c|c|}\hline
	 & {\tt Ordered}: & {\tt S2Lde} & {\tt L2Sde} & {\tt Ide} \\
	\hline
	Build 1 & 13.3 & 46.7 & 15.0 & 10.0 \\
    Build 2 & 21.7 & 48.3 & 16.7 & 15.0 \\
    \hline
\end{tabular}
\end{center}
\caption{Percentage of problems of dimension 4 and higher where each build type 
and heuristic resulted used the least number of iterations to 
detect a negative eigenvalue. (60 matrices after random shuffling of the 
variables).}
\label{tab:exactperm}
\end{table}

To further check how the use of these specific variables helps in structuring 
the Hessian matrix, we also applied a random orthogonal transformation to the 
variables in each problem prior to applying our heuristics.  The results 
appear in Table \ref{tab:exactorthog}, and concern 61 matrices\footnote{This 
number differs from the previous ones because applying an orthogonal transformation 
does not preserve the sign of the diagonal elements.}. In that setting, the
four strategies become more on par with one another using Build 1, while only the 
first two maintain a good performance using Build 2. Interestingly, 
combining Build 2 with the {\tt Ordered} heuristic emerges as the best variant.

\begin{table}[ht]
\begin{center}
\begin{tabular}{|l|c|c|c|c|c|}\hline
	 & {\tt Ordered}: & {\tt S2Lde} & {\tt L2Sde} & {\tt Ide} \\
	\hline
	Build 1 & 11.5 & 18.0 & 6.6 & 16.4 \\
	Build 2 & 49.2 & 41.0 & 9.8 & 16.4 \\
    \hline
\end{tabular}
\end{center}
\caption{Percentage of problems dimension 4 and higher where each build type 
and heuristic resulted used the least number of iterations to 
detect a negative eigenvalue (61 matrices after random orthogonal 
transformation).}
\label{tab:exactorthog}
\end{table}

\subsection{Using \NESAH{} on finite-difference Hessian approximations}
\label{ssec:numnesah}

In this section, we investigate the behavior of \NESAH{} applied with 
$\epsilon=0$. To this end, we repeat the experiment in 
Subsection~\ref{ssec:nesa} assuming we are only given access to the objective 
function of the 49 test problems. We thus compute finite-difference estimates 
of the 147 Hessian matrices from the previous section using the formulas given 
in Definition~\ref{def:DDHessian}, an oracle for the objective function as 
well as the points corresponding to the exact matrices described in 
Section~\ref{ssec:pbms}. For each matrix to approximate, we use three 
different values for the finite-difference parameter $h$, namely 
$\{ 10^{-2}, 10^{-4}, 10^{-6}\}$. We note that, due to approximation errors, 
the finite difference Hessian approximation failed to have negative curvature 
$4$ times for $h=10^{-2}$ and $3$ times for $h=10^{-6}$.
Overall, we obtain 231 test matrices with 
both negative curvature and no negative diagonal elements, 171 of which have 
dimension at least 4.  

Table~\ref{tab:NESAHall} presents our results when the 
four heuristics are used. Similarly to Tables~\ref{tab:exactALL} and 
\ref{tab:exact4up}, we observe that the {\tt Ordered} strategy leads to the 
fastest negative curvature detection. 

\begin{table}[ht]
    \centering
    \begin{tabular}{|l|l|c|c|c|c|c|}\hline
	Group &Matrices & {\tt Ordered} & {\tt S2Lde} & {\tt L2Sde} & {\tt Ide} \\\hline
    Build 1, Prob. Dim. $\ge 4$, all $h$ &171 &54.4 &38.0 &10.5 &14.6 \\
    Build 1, Prob. Dim. $\ge 4$, $h = 10^{-2}$ &56 &55.4 &37.5 &17.9 &19.6\\
    Build 1, Prob. Dim. $\ge 4$, $h = 10^{-4}$ &57 &54.4 &40.4 &3.5 &12.3\\ 
    Build 1, Prob. Dim. $\ge 4$, $h = 10^{-6}$ &58 &53.4 &36.2 &10.3 &12.1\\ \hline
    Build 2, Prob. Dim. $\ge 4$, all $h$ &171 &53.2 &47.4 &13.5 &23.4\\
    Build 2, Prob. Dim. $\ge 4$, $h = 10^{-2}$ &56 &53.6 &46.4 &16.1 &26.8 \\
    Build 2, Prob. Dim. $\ge 4$, $h = 10^{-4}$ &57 &52.6 &50.9 &10.5 &22.8\\  
    Build 2, Prob. Dim. $\ge 4$, $h = 10^{-6}$ &58 &53.4 &44.8 &13.8 &20.7\\ \hline
    \end{tabular}
    \caption{Percentage of problems of dimension 4 and higher where each build 
    type and heuristic resulted in the least number of iterations to detect a 
    negative eigenvalue (171 finite-difference matrices).}
    \label{tab:NESAHall}
\end{table}

We begin by commenting that the best variant of \NESAH{} detects negative curvature 
within $2$ iterations on $100$ matrices out of the $171$ considered in the tables, thereby 
using significantly less function evaluations than the upper bound provided in Lemma \ref{lem:NESAHterminates}. The worst case again corresponds to problem VAREIGVL
(dimension 10), where \NESAH{} uses $48$ function evaluations. This is still notably less 
than the $65$ function evaluations that would be required to approximate the entire Hessian.

As in the previous section, we repeat the experiment with a random reordering of 
the problem variables, so as to measure the significance of this ordering. 
Table~\ref{tab:permNESAHall} presents the results.  As in Section~\ref{ssec:nesa}, 
we observe that the performance of the {\tt Ordered} heuristic drops drastically 
in favor of the {\tt S2Lde} one.

\begin{table}[ht]
    \centering
    \begin{tabular}{|l|l|c|c|c|c|c|}\hline
	Group &Matrices & {\tt Ordered} & {\tt S2Lde} & {\tt L2Sde} & {\tt Ide} \\\hline
    Build 1, Prob. Dim. $\ge 4$, all $h$ &171 &12.9 &44.4 &12.3 &13.5 \\
    Build 1, Prob. Dim. $\ge 4$, $h = 10^{-2}$ &56 &16.1 &44.6 &19.6 &17.9\\
    Build 1, Prob. Dim. $\ge 4$, $h = 10^{-4}$ &57 &12.3 &45.6 &5.3 &10.5\\  
    Build 1, Prob. Dim. $\ge 4$, $h = 10^{-6}$ &58 &10.3 &43.1 &12.1 &12.1\\ \hline
    Build 2, Prob. Dim. $\ge 4$, all $h$ &171 &17.5 &46.2 &13.5 &20.5\\
    Build 2, Prob. Dim. $\ge 4$, $h = 10^{-2}$ &56 &19.6 &46.4 &16.1 &23.2 \\
    Build 2, Prob. Dim. $\ge 4$, $h = 10^{-4}$ &57 &17.5 &47.4 &10.5 &21.1 \\  
    Build 2, Prob. Dim. $\ge 4$, $h = 10^{-6}$ &58 &15.5 &44.8 &13.8 &17.2\\ \hline	
    \end{tabular}
    \caption{Percentage of problems of dimension 4 and higher where each build type 
	and heuristic resulted in the least number of iterations to 
	detect a negative eigenvalue (171 matrices after random shuffling of the 
	variables).}
    \label{tab:permNESAHall}
\end{table}

We also apply a random orthogonal transformation prior to applying our heuristics, 
leading to $235$ matrices with negative curvature and no negative diagonal elements. 
Among those matrices, $177$ are of dimension larger than or equal to $4$. The results 
appear in Table \ref{tab:orthogNESAHall}: once the orthogonal transformation is 
applied, we see all heuristics performing similarly for Build 1, whereas the 
combination of {\tt Ordered} and Build 2 outperforms the other variants.  
As such, we argue that the {\tt Ordered} 
heuristic is in fact the best choice of heuristic from this list.

\begin{table}[ht]
    \centering
    \begin{tabular}{|l|l|c|c|c|c|c|}\hline
	Group & Matrices &{\tt Ordered} & {\tt S2Lde} & {\tt L2Sde} & {\tt Ide} \\\hline
    Build 1, Prob. Dim.  $\ge 4$, all $h$ &177 &10.2 &16.4 &5.6 &16.9\\
    Build 1, Prob. Dim. $\ge 4$, $h = 10^{-2}$ &55 &10.9 &16.4 &5.5 &18.2\\
    Build 1, Prob. Dim. $\ge 4$, $h = 10^{-4}$ &61 &11.5 &18.0 &6.6 &16.4\\  
    Build 1, Prob. Dim. $\ge 4$, $h = 10^{-6}$ &61 &8.2 &14.8 &4.9 &16.4 \\ \hline
    Build 2, Prob. Dim. $\ge 4$, all $h$ &177 &50.8 &36.2 &11.3 &17.5\\
    Build 2, Prob. Dim. $\ge 4$, $h = 10^{-2}$ &55 &52.7 &34.5 &9.1 &18.2\\
    Build 2, Prob. Dim. $\ge 4$, $h = 10^{-4}$ &61 &50.8 &36.1 &13.1 &18.0\\  
    Build 2, Prob. Dim. $\ge 4$, $h = 10^{-6}$ &61 &49.2 &37.7 &11.5 &16.4 \\ \hline
    \end{tabular}
    \caption{Percentage of problems of dimension 4 and higher where each build type 
	and heuristic resulted in the least number of iterations to 
	detect a negative eigenvalue (177 finite-difference matrices with a random 
	orthogonal transformation).}
    \label{tab:orthogNESAHall}
\end{table}

\subsection{Comparison with all possible orderings}
\label{ssec:compall}

As a final experiment, we investigate whether there might exist another (as 
yet undiscovered) heuristic that might outperform the {\tt Ordered} heuristic. 
To do this, we ran \NESA{} and \NESAH{} with every possible ordering on 
matrices associated with the problems with dimensions between 4 and 8. Overall, 
we obtained 36 matrices without negative diagonal elements. 

Table~\ref{tab:allorders} details the results. For 26 problems out of 36, solving 
the problem using the {\tt Ordered} heuristics requires at worst two more iterations 
compared to using the {\tt Best} possible ordering. The {\tt Ordered} heuristic 
thus emerges as a very reasonable choice that is unlikely to be outperformed by any 
simple heuristic, especially in small dimensions.

\begin{table}[]
    \centering
    \begin{tabular}{|c|c|c|c|c|c|c|c|c|c|}
    \hline 
    Iterations for {\tt Ordered}    
    minus iterations for {\tt Best} &0 &1 &2 &3 &4 &5 &6 &7 &8\\
    \hline 
    Number of problems &11 &8 &7 &3 &3 &0 &2 &0 &2 \\
          \hline 
    \end{tabular}
    \caption{Comparison of the {\tt Ordered} heuristic with the best ordering 
    on 36 matrices (exact and finite-differences) with dimensions $4$ to $8$ 
    and no diagonal elements.}
    \label{tab:allorders}
\end{table}

\section{Concluding remarks}
\label{sec:conc}

We proposed an algorithm to detect negative eigenvalues in matrices, that 
proceeds by building submatrices given partial information.
The method is guaranteed to terminate in a finite number of steps, and comes 
with provable guarantees on its outputs. We have proposed eight variants of our 
algorithm based on two build strategies and four natural heuristics.  The variants were compared on a benchmark of 
Hessian matrices from the CUTEst library and their approximations through 
finite differences. Our experiments illustrate that considering the variables 
in order of appearance in the problem definition is often a insightful 
strategy to detect negative curvature without forming the entire matrix, that 
remains efficient on average when coupled with a strategy aiming at building 
large principle submatrices. Our software implementation, along with scripts 
to repeat the experiments herein, is available on github
\footnote{\href{https://github.com/clementwroyer/negative-eigs}
{https://github.com/clementwroyer/negative-eigs}}. The recommended settings 
apply Build 2 and the {\tt Ordered} heuristic.

Possible developments of our method would include dynamic modifications of our 
heuristics, that would exploit the new information available through the 
submatrices. Although our experiments suggest that the benefit might be limited in 
small dimensions, it could lead to provably lower error bounds.
Investigating the performance of our algorithm on approximate Hessians 
computed via alternate formulae is also an interesting avenue for future work.
Another perspective of this work consists in incorporating our algorithm 
into a DFO routine that could exploit negative curvature.  This will be the subject of future research.

\subsection*{Acknowledgements}

The authors are grateful to two anonymous referees, whose insightful comments 
lead to improvements in the \NESA{} and \NESAH{} algorithms.

\bibliographystyle{plain}
\bibliography{reffindiff}

\end{document}